\documentclass{amsart}

\usepackage{amsmath,amsthm,amssymb,amsfonts}

\newcommand{\beql}[1]{\begin{equation}\label{#1}}
\newcommand{\eeq}{\end{equation}}
\newcommand{\comment}[1]{}

\newcommand{\Abs}[1]{{\left|{#1}\right|}}

\newcommand{\Norm}[1]{{\left\|{#1}\right\|}}

\newcommand{\Set}[1]{{\left\{{#1}\right\}}}

\newcommand{\RR}{{\mathbb R}}

\newcommand{\ZZ}{{\mathbb Z}}

\newcommand{\dist}{{\rm dist\,}}

\newcommand{\ft}[1]{\widehat{#1}}

\newtheorem{theorem}{Theorem}

\newtheorem{lemma}{Lemma}

\theoremstyle{definition}
\newtheorem{remark}{Remark}
\newtheorem{example}{Example}

\newcounter{rem}
\newcommand{\rem}[1]{
       \refstepcounter{rem}
       \begin{quote} \noindent{\bf Comment \therem.}
       #1 \end{quote} }

\begin{document}

\title{A Weyl type formula for Fourier spectra and frames}

\author{Alex Iosevich \& Mihail N. Kolountzakis}
\date{October 2003}

\address{A.I.: Department of Mathematics, University of Missouri, Columbia, U.S.A}
\email{iosevich@math.missouri.edu}

\address{M.K.: Department of Mathematics, University of Crete, Knossos Ave.,
GR-714 09, Iraklio, Greece} \email{kolount@member.ams.org}

\thanks{The research of Alex Iosevich is partially supported by
the NSF Grant DMS02-45369. That of Mihalis Kolountzakis is
partially supported by European Commission IHP Network HARP
(Harmonic Analysis and Related Problems), Contract Number:
HPRN-CT-2001-00273 - HARP}

\begin{abstract}
We prove qualitative and quantitative results concerning the
asymptotic density in dilates of centered convex bodies of the
frequency vectors of orthogonal exponential bases and frames
associated to bounded domains in Euclidean space.
\end{abstract}

\maketitle

\section{Introduction}\label{sec:intro}
Let $\Omega \subset {\mathbb R}^d$ be a bounded domain. We say
that $\Omega$ is {\em spectral} if $L^2(\Omega)$ possesses an
orthogonal basis of the form $E_\Lambda={\{e^{2 \pi i x \cdot
\lambda}\}}_{\lambda \in \Lambda}$, where $\Lambda$ is a subset of
${\mathbb R}^d$. We shall refer to $\Lambda$ as a {\em spectrum}
of $\Omega$. The set $\Lambda$ is easily seen to have a separation
property: $\Abs{\lambda-\mu}\ge \epsilon$ for all
$\lambda,\mu\in\Lambda$, $\lambda\neq\mu$, and is therefore a
discrete set.

A systematic study of such spectra was initiated by Fuglede in
\cite{Fuglede74}. This problem has received much recent attention.

Let
$$
D_R^{+}(\Omega)=\max_{x \in {\mathbb R}^d} \# \{\Lambda \cap
Q_R(x)\},
$$
and
$$
D_R^{-}(\Omega)=\min_{x \in {\mathbb R}^d} \# \{\Lambda \cap
Q_R(x)\},
$$
where $Q_R(x)$ denotes the cube of side-length $2R$ centered at
$x$. Landau \cite{Landau67} proved that \beql{landau} \limsup_{R
\to \infty} \frac{D_R^{+}}{{(2R)}^d}=|\Omega|, \eeq and the same
equality holds for $D_R^{-}$.

If $E_\Lambda$ is only a frame for $L^2(\Omega)$, in the sense
that there exist positive constants $A$ and $B$ such that
\beql{frame} A{||g||}^2_{L^2(\Omega)} \leq \sum_A
{|\hat{g}(\lambda)|}^2 \leq B{||g||}^2_{L^2(\Omega)}, \eeq where
$\hat{g}$ denotes the Fourier Transform of $g$, then one can only
conclude that \beql{d-minus-frame} \limsup_{R \to \infty}
\frac{D_R^{-}}{{(2R)}^d} \ge |\Omega|. \eeq

Observe that $E_\Lambda$ is an orthogonal basis if and only if
\eqref{frame} holds with $A=B=1$.

In \cite{IoPe00}, the authors proved that if $\partial \Omega$ has
Minkowski dimension $\alpha<d$, with $\alpha$-dimensional upper
Minkowski content denoted by ${|\partial \Omega|}_{\alpha}$,
\eqref{frame} holds, and \beql{eq0.6} R \ge C {\left(
\frac{B{|\partial \Omega|}_{\alpha}}{A|\Omega|}
\right)}^{\frac{1}{d-\alpha}}, \eeq for some $C>0$, then for every
$x \in {\mathbb R}^d$, $Q_R(x)$ contains at least one element of
$\Lambda$.

In \S\ref{sec:cube} we show an alternative estimate
for the minimum $R$ such that every cube of
side-length $R$ necessarily contains a point of a $(A,B)$-frame
for the domain. For the quite general case of planar domains which
are simple polygons we show that this estimate is always as good
as \eqref{eq0.6} and often much better.

The main purpose of this paper is to prove a Weyl type estimate
(see e.g. \cite[Ch.\ 5]{Sogge93} for a description of the
classical Weyl asymtotics in the context of Riemannian manifolds)
for $\# \{\Lambda \cap R\cdot K\}$, where $K$ is a convex body in
${\mathbb R}^d$, symmetric with respect to the origin, in the case
when $E_\Lambda$ is an orthogonal basis for $L^2(\Omega)$.
In analogy with the
classical Weyl formula we ask whether
$$
\# \{\lambda: {||\lambda||}_K \leq
R\}=c_{\Omega,\Lambda}R^d+o(R^d),
$$
where ${||\cdot||}_K$ is the norm induced by $K$.

It turns out that the answer is yes, and under some additional
assumptions we obtain more quantitative estimates on the error
term. We shall see that our estimates imply both \eqref{landau}
and \eqref{eq0.6} at the same time, thus presenting a unified
description of geometric properties of spectra.

\begin{theorem}\label{th:main}
Let $\Omega$ be a bounded domain in ${\mathbb R}^d$ of positive
Lebesgue measure. Let $E_\Lambda$ be defined as above. If
$E_\Lambda$ is an orthogonal basis for $L^2(\Omega)$, then
\beql{eq0.8} \# \{\lambda \in \Lambda: {||\lambda||}_K \leq
R\}=|K|R^d|\Omega|+E(R), \eeq with \beql{qual} E(R)=o(R^d). \eeq
Moreover, if the boundary of $\Omega$ has upper Minkowski
dimension $\alpha<d$, then \beql{eq0.10} |E(R)| \leq C_K
{|\partial \Omega|}_{\alpha} R^{\alpha}, \eeq where $C$ is a
constant depending on the domain $\Omega$ only.

If $E_\Lambda$ is merely a frame with constants $A$ and $B$ as in
\eqref{frame} above, then \beql{frame-conclusion} A|\Omega| \leq
\frac{\# \{\lambda \in \Lambda: {||\lambda||}_K \leq R \}}{|K|R^d}
\leq B|\Omega|. \eeq
\end{theorem}

\begin{remark} We do not have an example where the estimate
$|E(R)| \lesssim R^{\alpha}$ holds only with $\alpha>{d-1}$. A
family of examples with $\alpha=d-1$ is given in the Example 1
below. Nevertheless, the upper Minkowski assumption allows to
obtain better results than the $o(R^d)$ estimate obtained without
any assumption on the boundary of $\Omega$. \end{remark}

Observe that \eqref{eq0.6} follows immediately from \eqref{eq0.8}
and  \eqref{eq0.10} in the case $\Lambda$ is a spectrum, while
\eqref{landau} and \eqref{d-minus-frame} follow from \eqref{qual}
and \eqref{frame-conclusion}, respectively. To see the former, we
use the fact that the proof of \eqref{eq0.8} and \eqref{eq0.10} shows
that the same estimate still holds with uniform constants for a
translated spectrum, i.e $\# \{\lambda \in \Lambda:
{||\lambda-x||}_K \leq R\}=|K|R^d|\Omega|+E(x,R)$, with $|E(x,R)|
\leq C_K {|\partial \Omega|}_{\alpha} R^{\alpha}$ with constants
independent of $x$, under the assumption that the boundary of
$\Omega$ has upper Minkowski dimension $\alpha<d$. It follows that
if $C_K {|\partial \Omega|}_{\alpha} R^{\alpha}$ is much smaller
than the main term $|K|R^d|\Omega|$, the set $\{\lambda \in
\Lambda: {||\lambda-x||}_K \leq R\}$ is not empty and
\eqref{eq0.6} follows.

It is interesting to compare Theorem \ref{th:main} with a model
situation where $\Omega={[0,1]}^d$, $K$ is the unit Euclidean ball
of volume $\omega_d$, and $\Lambda={\mathbb Z}^d$. In this case,
it is known that \beql{eq0.12} \# \{\lambda \in \Lambda: |\lambda|
\leq R\}=\omega_d R^d +O(R^{d-2}), \eeq if $d \ge 5$, and
\beql{eq0.13} \# \{\lambda \in \Lambda: |\lambda| \leq
R\}=\omega_d R^4 +O(R^{2}\log(R)), \eeq if $d=4$. In two and three
dimensions the situation is more murky. In two dimensions, Hardy's
conjecture says that the remainder should be
$O(R^{\frac{1}{2}+\epsilon})$, for any $\epsilon>0$. The best
known result in dimension two, due to Huxley, gives
$O(R^{\frac{46}{73}})$. In three dimensions, the best known result
is due to Heath-Brown who obtained $O(R^{\frac{21}{16}})$. See
\cite{Huxley96} and the references contained therein.

All these results are driven by the curvature of the boundary of
the Euclidean ball. To see this observe that if
$\Omega={[0,1]}^d$, $\Lambda={\mathbb Z}^d$, and $K={[-1,1]}^d$,
the best we can say is \beql{eq0.14} \# \{\lambda \in \Lambda:
{||\lambda||}_K \leq R\}=2^d R^d +O(R^{d-1}), \eeq since an
integer dilate of the unit cube contains $\approx R^{d-1}$ integer
lattice points on its boundary.

Observe that the remainders in all these results are much better
than the ones given by Theorem \ref{th:main}, especially in the
case where $K$ is the Euclidean ball. However, the following
construction shows that even when $\Omega$ is the unit cube in
${\mathbb R}^d$ and $K$ is the Euclidean ball, the situation
becomes much worse when the integer lattice is replaced by a more
complicated spectrum.

\begin{example}
Consider the standard lattice of cubes in ${\mathbb R}^d$. View
each cube as a translate of ${[0,1]}^d$. We call the image of the
origin under this translation the defining vertex of the cube.

Consider a sphere of radius $100$ centered at the origin. Move
each relevant column of cubes in such a way that the defining
vertex of one the cubes in the column is on this sphere. Now
consider a sphere radius $2^{100}$. Leave the previously moved
columns alone, and move the other columns in such a way that the
defining vertex of one of the cubes in each column is on this
sphere. Continuing this process, we produces a family of spheres
with radii $\{R_i\}$, $R_i \to \infty$, and a discrete set
$\Lambda'$ such that the $i$'th sphere contains $\approx
R_i^{d-1}$ points of $\Lambda'$ with constants independent of $i$.

Since ${[0,1]}^d+\Lambda'$ is clearly a tiling, the result proved
independently by Lagarias, Reeds and Wang \cite{LRW98} and
Iosevich and Pedersen \cite{IoPe98} (see also \cite{k}), implies
that $E_{\Lambda'}$ is an orthogonal basis for $L^2(\Omega)$. It
is clear that in this case, the estimate \eqref{eq0.10} cannot be
improved.
\end{example}

\section{Proof of Theorem \ref{th:main}}\label{sec:main-proof}
Suppose that $E_\Lambda$ is an orthogonal basis for $L^2(\Omega)$.
It follows (see e.g.\ \cite{k} or Remark \ref{rem:frames} below) that
\beql{tiling}
\sum_{\lambda\in\Lambda} f(x-\lambda)= \Abs{\Omega}^2, \eeq where
$f = \Abs{\ft{\chi_\Omega}}^2$. Notice also that $\int f =
\Abs{\Omega}$.

Let $K_t$ denote the set $t\cdot K$, let $\Lambda_t$ denote the
set $\Lambda \cap K_t$ and let finally
$$
N(t) = \#\Lambda_t.
$$
From the orthogonality of the exponentials in $E_\Lambda$ it
follows that for any two distinct points $\lambda,\mu\in\Lambda$
their difference is a zero of $f$, hence there is a lower bound
for their distance. This separation means that $N(t) \le C t^d$
and also that $N(t+R)-N(t) \le C t^{d-1} R$ if $R>0$ is a large
enough constant. Roughly speaking, if a domain is not too thin
then the number of $\Lambda$-points in it is bounded above by its
volume.

\begin{remark}\label{rem:frames}
If $E_\Lambda$ is only a frame for the space $L^2(\Omega)$ with
frame constants $A$ and $B$ then it follows that, for
$f=\ft{\chi_\Omega}$ we have for almost all $x \in \RR^d$
\beql{packing-covering} A\Abs{\Omega}^2 \le
\sum_{\lambda\in\Lambda} f(x-\lambda) \le B\Abs{\Omega}^2. \eeq
This is easily seen by applying \eqref{frame} to an arbitrary
exponential function. The proof below for the case of $E_\Lambda$
being an orthonormal basis also gives \eqref{frame-conclusion}
with trivial modifications: the separation property does not hold
for frames, however using the upper bound in
\eqref{packing-covering} one easily gets that $\Lambda$ has the
property that matters in the proof below, namely that the number
of $\Lambda$-points in $K$-balls and shells is controlled by the
volume of the region.
\end{remark}

First we follow \cite{kl} to show that $N(t) = \Abs{K}
\Abs{\Omega} t^d + o(t^d)$, under no assumptions about
$\partial\Omega$. For this let $1<R<T$ be large numbers and
integrate \eqref{tiling} over the region $K_T$ to get
\begin{eqnarray*}
\Abs{\Omega}^2 \Abs{K_T}
&=& \int_{K_T} \sum_{\lambda\in\Lambda} f(x-\lambda)~dx\\
&=& \int_{K_T}\sum_{\Lambda_{T-R}} f(x-\lambda)~dx +
    \int_{K_T}\sum_{\Lambda_{T+R}\setminus\Lambda_{T-R}} f(x-\lambda)~dx +
    \int_{K_T}\sum_{\Lambda\setminus\Lambda_{T+R}} f(x-\lambda)~dx \\
&=& \int_{\RR^d}\sum_{\Lambda_{T-R}} f(x-\lambda)~dx - E_1 + E_2 + E_3\\
&=& N(T-R)\Abs{\Omega} - E_1 + E_2 + E_3,
\end{eqnarray*}
so that
$$
N(T-R) = \Abs{K_T}\Abs{\Omega} -\frac{E_1}{\Abs{\Omega}}+
               \frac{E_2}{\Abs{\Omega}}+
               \frac{E_3}{\Abs{\Omega}},
$$
where
\begin{eqnarray*}
E_1 &=& \int_{K_T^c} \sum_{\Lambda_{T-R}} f(x-\lambda)~dx,\\
E_2 &=& \int_{K_T} \sum_{\Lambda_{T+R}\setminus\Lambda_{T-R}} f(x-\lambda)~dx,\\
E_3 &=&
\int_{K_T}\sum_{\Lambda\setminus\Lambda_{T+R}}f(x-\lambda)~dx.
\end{eqnarray*}
Let $\epsilon>0$ be arbitrary and fixed, and choose $R$ so that
$$
\int_{K_R^c} f(x)~dx \le \epsilon.
$$
We have, since $T>R$,
\begin{eqnarray*}
E_1 &\le& \epsilon N(T-R)\\
&\le& \epsilon N(T)\\
&\le& C \epsilon T^d,
\end{eqnarray*}
and
\begin{eqnarray*}
E_3 &=&\sum_{\Lambda\setminus\Lambda_{T+R}} \int_{\lambda+K_T} f(x)~dx\\
 &\le& C T^d \int_{K_R^c} f(x)~dx \\
 &\le& C \epsilon T^d,
\end{eqnarray*}
as each point in $K_R^c$ is contained in at most $CT^d$ of the
sets $\lambda+K_T$, with
$\lambda\in\Lambda\setminus\Lambda_{T+R}$. For $E_2$ we trivially
have
\begin{eqnarray*}
E_2 &\le& (N(T+R) - N(T-R)) \Abs{\Omega}\\
   &=& o(T^d),\ \ \mbox{as $R$ is fixed}.
\end{eqnarray*}
Since $\epsilon$ is arbitrary and $N(T-R)=N(T)-o(T^d)$ as $T \to
\infty$ we have proved that
$$
N(T) = \Abs{K}\Abs{\Omega}T^d + o(T^d).
$$
In other words the set $\Lambda$ has density $\Abs{\Omega}$.

Let us now assume that $\partial\Omega$ has finite
$\alpha$-dimensional upper Minkowski content, with $\alpha<d$. We
shall prove that \beql{decay} \int_{K_{2R}\setminus K_R} f(x)~dx
\le
 C_K \Abs{\partial\Omega}_\alpha R^{-(d-\alpha)},
\eeq where $\Abs{\partial\Omega}_\alpha$ denotes the Minkowski
content of $\partial\Omega$.

To see \eqref{decay} we choose $N$ (independent of $j$) boxes $Q_\nu$ and $N$
vectors $h_\nu$, with $2^{-(j+1)}\le \Norm{h_\nu}_K \le 2^{-j}$, such
that the region $\Set{y:\ 2^j \le \Norm{y}_K \le 2^{j+1}}$ is
contained in $\bigcup Q_\nu$ and
$$
\Abs{e^{2\pi i y\cdot h_\nu} -1} \ge 1,\ \ (y\in Q_\nu).
$$
We now have, where $2^j \sim R$,
\begin{eqnarray*}
\int_{2^j \le \Norm{y}_K \le 2^{j+1}} f(y)~dy
& \le & \sum_{\nu=1}^N \int_{Q_\nu}
          \Abs{\ft{\chi_\Omega}(y)(e^{2\pi i y\cdot h_\nu} -1)}^2 ~dy\\
& \le & \sum_{\nu=1}^N \int_{\RR^d}
      \Abs{\chi_\Omega(x) - \chi_\Omega(x-h_\nu)}^2~dx\\
& \le & \sum_{\nu=1}^N \Abs{\Set{x\in\RR^d:\ \dist(x,\partial\Omega) < h_\nu}}\\
& \le & C \sum_{\nu=1}^N \Abs{\partial\Omega}_\alpha \Abs{h_\nu}^{d-\alpha}\\
& \le & C_K \Abs{\partial\Omega}_\alpha 2^{-j(d-\alpha)}.
\end{eqnarray*}

Let $R\to\infty$, $\epsilon=\int_{K_R^c}f=O(R^{-(d-\alpha)})$ and
integrate \eqref{tiling} on $K_R$ to get, similarly to what we did
above,
$$
N(R) = \Abs{\Omega} \Abs{K_R} -\frac{E_1}{\Abs{\Omega}}
                 +\frac{E_2}{\Abs{\Omega}}
                 +\frac{E_3}{\Abs{\Omega}},
$$
where now
\begin{eqnarray*}
E_1 &=& \int_{K_R^c} \sum_{\Lambda_{R}} f(x-\lambda)~dx,\\
E_2 &=& \int_{K_R} \sum_{\Lambda_{2R}\setminus\Lambda_{R}} f(x-\lambda)~dx,\\
E_3 &=&
\int_{K_R}\sum_{\Lambda\setminus\Lambda_{2R}}f(x-\lambda)~dx.
\end{eqnarray*}
We have as before
$$
E_3 \le R^d \int_{K_R^c} f = O(R^\alpha).
$$
To bound $E_1$ we decompose the set $\Lambda_R$ in shells of width
$2^j$, $j=0,1,\ldots,\log_2 R$:
$$
A_j = \Lambda \cap (K_{R-2^j} \setminus K_{R-2^{j+1}}),
$$
thinner near $\partial K_R$ and doubling in width as we move
towards the origin. Using \eqref{decay} we get
\begin{eqnarray*}
E_1 &\le& \sum_{j=0}^{\log_2 R} \# A_j 2^{-j(d-\alpha)}\\
&\le& C \sum_{j=0}^{\log_2 R} R^{d-1} 2^j 2^{-j(d-\alpha)}\\
&\le& C R^{d-1} \sum_{j=0}^{\log_2 R} 2^{j(1-d+\alpha)}\\
&=& O(R^\alpha).
\end{eqnarray*}
We bound $E_2$ similarly by decomposing in dyadic shells the set
$\Lambda_{2R}\setminus \Lambda_R$, thinner near $\partial K_R$ and
doubling as we move out:
$$
B_j = \Lambda \cap (K_{R+2^{j+1}} \setminus K_{R+2^j}),
$$
We get
\begin{eqnarray*}
E_2 &\le& \sum_{j=0}^{\log_2 R} \# B_j 2^{-j(d-\alpha)}\\
&\le& C \sum_{j=0}^{\log_2 R} (R+2^j)^{d-1} 2^j 2^{-j(d-\alpha)}\\
&\le& C R^{d-1} \sum_{j=0}^{\log_2 R} 2^{j(1-d+\alpha)}\\
&=& O(R^\alpha).
\end{eqnarray*}
We have proved
$$
N(R) = \Abs{\Omega}\Abs{K}R^d + O(R^\alpha).
$$

\section{An improved upper bound for the side-length of empty cubes}
\label{sec:cube}

In this section we show a very simple new upper bound on the size
of $\Lambda$-free cubes, when $E_\Lambda$ is a frame for
$L^2(\Omega)$. This bound ignores the roughness of
$\partial\Omega$ and cares about the thickness of the interior.
\begin{theorem}\label{th:trivial}
Suppose $\Omega\subseteq\RR^d$ is measurable with positive measure
and also that the point set $\Lambda\subseteq\RR^d$ is such that
$E_\Lambda$ is a frame for $L^2(\Omega)$ with lower and upper
frame constants being $A$ and $B$, respectively. Let $R$ be such
that there is a cube in $\RR^d$ of side-length $R$ containing no
point of $\Lambda$. Finally, assume that $\Omega$ contains a cube
of side-length $\epsilon>0$. Then
\begin{equation}\label{new-estimate}
R \le C_d \frac{B}{A} \epsilon^{-1}.
\end{equation}
\end{theorem}
\begin{proof}
Let $Q \subseteq \Omega$ be a cube of side-length $\epsilon$
contained in $\Omega$. By the definition of a frame \eqref{frame}
it follows that $E_\Lambda$ is a frame for
$L^2(Q)$ with the same frame constants. Applying \eqref{eq0.6} for
$Q$, with $\alpha=d-1$, we get the result.
\end{proof}

The new estimate \eqref{new-estimate} can easily seen to be much
better than \eqref{eq0.6} in some cases. Take for example in
$\RR^2$ a perturbation of the unit square that creates a very long
boundary but leaves intact a square of side $1/2$ inside. And, in
the quite general case of simple polygonal domains the new
estimate \eqref{new-estimate} is at least as good as
\eqref{eq0.6}, as the following theorem claims.

\begin{theorem}\label{th:comparison}
Suppose $\Omega$ is a simple polygon and let $\epsilon$ be the
maximum side-length of a cube contained in $\Omega$. Then, for
some constant $C>0$, we have
\begin{equation}\label{comp}
\frac{\Abs{\partial\Omega}}{\Abs{\Omega}} \ge C \epsilon^{-1}.
\end{equation}
\end{theorem}

\begin{proof}
Consider the usual subdivision of the plane into squares of side
$\epsilon$, translated at the points $(\epsilon\ZZ)^2$. Let $N$ be
the number of those squares intersecting $\Omega$. It follows that
$N \ge \Abs{\Omega}\epsilon^{-2}$.

By our assumption about $\epsilon$ these squares all contain some
point of $\Omega^c$ and therefore also some point of
$\partial\Omega$. Let these squares be called $Q_i$,
$i=1,\ldots,N$, and let $p_i \in Q_i \cap \partial\Omega$.

Partition the set of these squares into four classes depending on
the parity of the $x$- and $y$-coordinate of their lower left
corner (after multiplying these coordinates by $\epsilon^{-1}$).
At least one of these classes contains at least $N/4$ squares.
There are therefore at least $N/4$ points on $\partial\Omega$ with
minimum distance $C \epsilon$ from each other, and these are
connected to each other along $\partial\Omega$. It follows that
the total length of $\partial\Omega$ is at least
$$
C \frac{N}{4} \epsilon \ge C \Abs{\Omega} \epsilon^{-1}.
$$
\end{proof}

\enddocument